\documentclass{amsart}
\usepackage{amssymb, amsfonts} 
\usepackage[v2]{xy}

\newtheorem{thm}{Theorem}[section]
\newtheorem{cor}[thm]{Corollary}
\newtheorem{prop}[thm]{Proposition}
\newtheorem{lem}[thm]{Lemma}

\theoremstyle{remark}
\newtheorem{warn}[thm]{Warning}

\theoremstyle{definition}
\newtheorem{defn}[thm]{Definition}

\newcommand{\sT}{\mathcal T}

\newcommand{\J}{\mathcal C}

\newcommand{\bQ}{\mathbb{Q}_{X}}
\newcommand{\spec}{Sp^{\Sigma}}
\newcommand{\overto}[1]{\xrightarrow{#1}}
\newcommand{\rtarr}{\longrightarrow}
\newcommand{\com}{\circ}
\newcommand{\ul}{\underline}
\newcommand{\sma}{\wedge}
\newcommand{\htp}{\simeq}
\newcommand{\si}{\sigma}
\newcommand{\al}{\alpha}
\newcommand{\ph}{\phi}
\newcommand{\ze}{\zeta}
\newcommand{\et}{\eta}
\newcommand{\epz}{\varepsilon}
\newcommand{\ta}{\tau}
\newcommand{\SI}{\Sigma}
\DeclareMathOperator{\colim}{colim}
\DeclareMathOperator{\map}{map}
\DeclareMathOperator{\tel}{tel}
\DeclareMathOperator{\id}{id}
\DeclareMathOperator{\Id}{Id}
\DeclareMathOperator{\hocolim}{hocolim}
\newcommand{\iso}{\cong}
\newcommand{\OM}{\Omega}
\newcommand{\la}{\lambda}

\makeatletter
\let\c@equation\c@thm
\makeatother
\numberwithin{equation}{section}

\begin{document}

\title{A telescope comparison lemma for THH}

\author{M.A. Mandell}
\thanks{Both authors partially supported by NSF Postdoctoral Fellowships}
\address{Department of Mathematics\\ MIT \\ 77  Massachusetts Ave. \\ Cambridge, MA 02139}
\email{mandell@math.mit.edu}
\author{B. Shipley}
\address{Department of Mathematics \\ University of Chicago \\ 
5734 S. University Ave. \\ Chicago, IL 60637}
\email{bshipley@math.uchicago.edu}

\begin{abstract}
We extend to the non-connective case a lemma of B\"okstedt about the
equivalence of the telescope with a more complicated homotopy colimit
of symmetric spectra used in the construction of THH.
\end{abstract}
\maketitle

\bibliographystyle{plain}

\maketitle

\section{Introduction}

In \cite{bo}, B\"okstedt constructed the topological Hochschild homology
of a symmetric ring spectrum, before the smash product of symmetric spectra
had been invented.  Now that we understand the smash product of
symmetric spectra, we can construct THH using the cyclic bar
construction.  The paper \cite{ship} compares these definitions of THH.
Both of these constructions are the geometric realizations of
simplicial symmetric spectra that turn out to be stably equivalent in
each simplicial level.  A peculiarity of the theory of symmetric
spectra is that not all stable equivalences are equivalences of the
underlying prespectra, and simple examples show that these two
constructions of THH do not always have weakly equivalent underlying
prespectra.

The problem can be understood levelwise, where the ring structure
plays no role.  (The multiplication is only used to define the
simplicial face maps.)  Understanding simplicial level zero is the key to
understanding all the levels, and we concentrate on this.  As
traditional in working with THH, instead of thinking in terms of
prespectra, we can think in terms of functors from based spaces to
based spaces.  Then for a 
symmetric spectrum $T$, the zeroth level of the cyclic bar construction
$thh_{0}(T)$ is just $T$ itself, and this corresponds to the functor
that sends a based space $X$ to the space $\tel_{n}
\Omega^{n}(T_{n}\sma X)$, that is, this space is the underlying
infinite loop space of $T\sma X$.  The zeroth level
of B\"okstedt's construction is a homotopy colimit.  It sends $X$ to
the space 
$$ THH_{0}(T;X)=\hocolim_{n\in I} \Omega^{n}(T_{n}\sma X). $$
Here $I$ is the category finite sets and injections; see
Definition~\ref{def I}.  There is a canonical inclusion of the
telescope in this homotopy colimit and B\"okstedt's fundamental lemma
for THH is that it is a weak equivalence when $T$ is strictly
connective ($\pi_{q}T_{n}=0$ for $q<n$) and convergent
(Definition~\ref{dd1}).  The purpose of this paper is to eliminate the
connectivity hypothesis. 

This work grew out of the authors' previous study of topological Hochschild
homology as defined for symmetric spectra and for $S$-modules in~\cite{ekmm}
and~\cite{ship}.  This material also formed the basis of a first comparison
of the homotopy categories of symmetric spectra, $S$-modules, and several
other new symmetric monoidal categories of spectra carried out in a larger
project with May and Schwede.  See~\cite{may} for an expository account
based on this earlier outline.  This first approach has since been
replaced by highly structured comparisons of model categories,
see~\cite{mmss} and~\cite{Sch}.

The definition of symmetric spectra easily generalizes to other, even
non-top\-o\-lo\-gi\-cal\-ly motivated contexts.  For example,
symmetric spectra have also been recently used to model stable
${\mathbb A}^1$-local homotopy theory in the sense
of~\cite{voevodsky-ICM}, q.v.~\cite{jardine}. Although the proofs do
not immediately transport to this setting, one expects similar
statements to hold.  In fact, the impetus for resurrecting this
material came from questions posed by Voevodsky that are answered
here.

We work in the setting of topological symmetric spectra, but
only minor modifications, such as level fibrant replacement, would be
necessary to apply these results to simplicial symmetric spectra.  To
avoid topological technicalities, we assume that any level of a given
symmetric spectrum or any given based space here is non-degenerately based.

\subsection*{Acknowledgments:} As mentioned above, this material was
developed in the context of a larger project with Peter May and Stefan
Schwede.  We would like to thank them for their support during this
work and for encouraging us to publish this piece of work separately.
In particular, we would like to thank Peter May for his editorial and
mathematical guidance throughout the writing of this paper.

\section{Statement of results}
We first need to recall a few definitions before we state our results.

\begin{defn}\label{dd1}
Let $T$ be a symmetric spectrum. It has an underlying prespectrum with
spaces $T_n$ and structure maps $\si_n: T_n\sma S^1 \rtarr T_{n+1}$.
Define the homotopy groups of $T$ to be the homotopy groups of the prespectrum:
$$\pi_qT = \colim \pi_{q+n}T_n \iso \pi_q(\tel_n\OM^n T_n).$$
We say that $T$ is {\em convergent} if the adjoint structure map 
$\tilde{\si}_n: T_n\rtarr \OM T_{n+1}$ is an $(n+\la(n))$-equivalence for 
some nondecreasing sequence of integers $\la(n)$ with infinite limit. We say that 
$T$ is a {\em symmetric $\OM$-spectrum} if the maps $\tilde{\si}_n$ are weak 
equivalences. We say that a map $f:T\rtarr T'$ of symmetric spectra is a 
{\em $\pi_*$-isomorphism} if it induces an isomorphism of homotopy groups. 
We say that $f$ is a {\em level weak equivalence} if each 
$f_n: T_n \rtarr T'_n$ is a weak equivalence. 
Let $[X,Y]$ denote the set of maps $X\rtarr Y$ in the homotopy 
category with respect to the level model structure \cite[\S6]{mmss}.
We say that $f$ is a  {\em stable weak equivalence} if 
$f^* : [Y,E] \rtarr [X,E] $ 
is a bijection for all symmetric $\OM$-spectra $E$. 
\end{defn}

Recall from \cite{hss, mmss} that a $\pi_*$-isomorphism is a stable weak equivalence,
but not conversely. However, the converse does hold for stable weak equivalences between 
convergent spectra, and we shall only be concerned with convergent spectra here.

\begin{defn}\label{def I} Let $I$ be the category whose objects are the sets 
$\mathbf n =\{1,\ldots, n\}$, $n\geq 0$, and whose morphisms are 
the injective functions. Let $\SI$ be the subcategory
of isomorphisms in $I$. For $m < n$, let $i_{m,n}: \mathbf m \rtarr \mathbf n$
be the standard inclusion and observe that any map $f: \mathbf m\rtarr \mathbf n$ 
can be written, non-uniquely, as a composite $f=\ta \circ i_{m,n}$ for some 
$\ta \in \SI_{n}$. Let $J\subset I$ be the subcategory of objects $\mathbf n$ 
and standard injections $i_{m,n}: \mathbf m\rtarr \mathbf n$ for $m\leq n$.
\end{defn}

Our main theorem is the following sharpening of a basic lemma of B\"okstedt \cite[1.6]{bo}
that eliminates its connectivity hypotheses; see also \cite[3.1.7]{ship} for a further generalization in the case $X=S^0$. Let $\sT$ be the
category of based spaces. A sequence of based spaces $X_n$ and maps $X_n\rtarr X_{n+1}$ 
determines a functor $J\rtarr \sT$, and there is a natural homotopy equivalence 
$$\tel_n X_n \overto{\htp} \hocolim_J X_n.$$
A functor defined on $I$ restricts to a functor defined on $J$ and thus 
gives an induced map of homotopy colimits. 

A symmetric spectrum $T$ gives rise to a functor $I\rtarr \sT$ that sends $\mathbf n$ to 
$\OM^{n}(T_n\sma X)$. Using the evident natural map 
$\OM^n(X)\sma Y \rtarr \OM^n(X\sma Y)$, the adjoint structure map $\tilde\si$ gives
$$T_m\sma X\rtarr (\OM^{n-m}T_n)\sma X\rtarr \OM^{n-m}(T_n\sma X);$$
the map induced by $i_{m,n}$ is obtained from this by applying 
$\OM^{m}$. Permutations act on $S^n$ and $T_n$ and by conjugation on 
$\OM^n(T_n\sma X)$. This construction plays a fundamental role in the 
theory of symmetric spectra, see~\cite{ship}. 

\begin{thm}[B\"okstedt's lemma]\label{boklem}
Let $T$ be a convergent symmetric spectrum and $X$ be a based
space. Then the natural map 
$$ \tel_n \OM^{n}(T_n\sma X)\htp \hocolim_J\OM^{n}(T_n\sma X)
\rtarr \hocolim_I\OM^{n}(T_n\sma X) $$
is a weak equivalence.
\end{thm}

A proof of B\"okstedt's original lemma appears in \cite[2.3.7]{Madsen}.
When the symmetric spectrum $T$ is not assumed to be connective, the
homotopy groups of $T_{n}\sma X$ do not stabilize, and the proof of
Theorem~\ref{boklem} requires the construction and analysis of
$(-k)$-connected covers; see Lemma~\ref{key}.  B\"okstedt and Madsen's
argument applies on each cover, and we use standard techniques
with homotopy colimits to mesh them together and deduce the result for
$T$.  For the reader's convenience we review the definition and basic
properties of homotopy colimits in Section \ref{HOCO}. 

It is a standard fact in stable homotopy theory that the map
$$ \tel_n \OM^{n}(T_n\sma X)\to \OM(\tel_{n} \OM^{n}(T_n\sma \SI X))  $$
is a weak equivalence; see for example \cite[7.4.i']{mmss}.  It
follows that when $T$ is convergent, the map
$$ \hocolim_{I} \OM^{n}(T_n\sma X)\to \OM(\hocolim_{I} \OM^{n}(T_n\sma \SI X))  $$
is a weak equivalence.  We obtain the following corollary of
Theorem~\ref{boklem}. 

\begin{cor}\label{zeroth} If $T$ is convergent, then, for a based
space $X$, we have an isomorphism of homology 
$$ T_{q}(X)=\pi_{q}(T\sma X) \iso \pi_{q+k} \hocolim_{I}
\Omega^{n}(T_{n}\sma \SI^{k}X) $$
and an isomorphism of cohomology 
$$ T^{q}(X)=\pi_{-q}(F(X,T))\iso \pi_{k-q}F(\SI^{k-q} X, \hocolim_{I}
\Omega^{n}(T_{n}\sma S^{k})). $$
\end{cor}

For $S$ the sphere spectrum, $S\sma X$ is isomorphic to $\SI^{\infty}X$, 
and we obtain the following corollary. 
Recall that $QX = \colim \OM^n\SI^n X$. 
Since the maps of the colimit system are inclusions, the natural map 
$\tel_n \OM^n\SI^n X\rtarr QX$ is a weak equivalence. 

\begin{cor}\label{Susp} The natural map
$$\tel_n \OM^{n}\SI^n X\htp \hocolim_J\OM^{n}\SI^n X
\rtarr \hocolim_I\OM^{n}\SI^n X $$
is a weak equivalence.
\end{cor}

\section{The proof of B\"okstedt's lemma}

The strategy for the proof of Theorem \ref{boklem} is to reduce to the
study of a symmetric $\OM$-spectrum that is closely related to $T\sma X$. 
Let $\bQ T$ be the symmetric spectrum with levels  
$$ \bQ T_n=\tel_{m} \OM^{m} (T_{n + m}\sma X). $$
Then $\bQ$ is a functor of both $T$ and $X$.  The inclusions of the 
initial term in a telescope gives a natural map of symmetric spectra 
$T\sma X\rtarr \bQ T$ such that the following diagram commutes:
$$ \diagram
\tel \OM^{n}(T_n\sma X)\rto\dto
&\tel \OM^{n}\bQ T_n \dto\\
\hocolim_I \OM^{n}(T_n\sma X)\rto
&\hocolim_I \OM^{n}\bQ T_n.
\enddiagram $$

Theorem~\ref{boklem} states that the left vertical
arrow is a weak equivalence. To prove this, we prove 
that the remaining three arrows are weak equivalences.  We prove the
following two lemmas in the next section.

\begin{lem}\label{lemma1}
If $T$ is a convergent symmetric spectrum, then $\bQ T$ is a symmetric
$\OM$-spectrum and the natural map $T\sma X\rtarr \bQ T$ is a
$\pi_*$-isomorphism.  
\end{lem}

For $q\geq 0$, the map of $q$th homotopy groups of spaces induced by
the top arrow in the diagram is the map of $q$th homotopy groups of
prespectra induced by $T\sma X\rtarr \bQ T$, so this lemma implies
directly that the top arrow is a weak equivalence.

\begin{warn}\label{warn}
If $T$ is not convergent, then $\bQ T$ is not a symmetric $\OM$-spectrum 
and $T\sma X\rtarr \bQ T$ is not a $\pi_{*}$-isomorphism in general, even 
when $X=S^{0}$. In fact, it is shown in \cite{hss} that the symmetric spectrum 
$T=F_1(S^1)$ defined below gives a counterexample.
\end{warn}

The next lemma shows that the bottom
horizontal arrow in the diagram is a  weak equivalence.

\begin{lem}\label{lemma2}
If $T$ is a convergent symmetric spectrum, then the natural map
$$ \hocolim_I \OM^{n}(T_n\sma X)\rtarr
\hocolim_I\OM^{n}\bQ T_n $$
is a weak equivalence.
\end{lem}

Finally, the following lemma 
allows us to exploit the fact that $\bQ T$ 
is a symmetric  $\OM$-spectrum to prove that the right vertical arrow in 
the diagram is a weak equivalence.

\begin{lem}\label{lemma3}
If $U$ is a symmetric  $\OM$-spectrum, then the natural map
$$U_0\rtarr \hocolim_{I}\OM^{n}U_n $$
is a weak equivalence. 
\end{lem}

This lemma follows from a standard result about homotopy colimits of
homotopically constant diagrams, see Lemma \ref{includezero}. 

By Lemmas \ref{lemma1} and \ref{lemma3}, when $T$ is convergent the inclusion 
$$\bQ T_0\rtarr \hocolim_I\OM^{n}\bQ T_n$$ 
is a weak equivalence. The inclusion 
$$\bQ T_0\rtarr \tel \OM^{n}\bQ T_n$$ 
is also a weak equivalence by Lemma~\ref{includezero}. Since the right 
vertical arrow is a map under
$\bQ T_0$, it too is a weak equivalence. This completes the proof
of Theorem \ref{boklem}, modulo the proofs of the lemmas.

\section{The proofs of two technical lemmas}

The proofs of Lemmas \ref{lemma1} and \ref{lemma2} are based on the use of
``strictly $k$-connected covers'' of symmetric spectra.

\begin{defn}  
A symmetric spectrum $T$ is {\em strictly $k$-connected} if
$T_n$ is $(n+k)$-connected for $n\geq 0$.
\end{defn}

We prove the following lemma in the next section. 

\begin{lem}\label{covering}
For each integer $k$, there is a functor $C_{k}$ on the category of
symmetric spectra such that, for any symmetric spectrum $T$, $C_{k}T$ is
strictly $k$-connected.  There are natural transformations
$c_{k}: C_{k}\rtarr\Id$ and $e_{k}: C_{k}\rtarr C_{k-1}$ such
that $c_{k-1}=c_{k}\circ e_{k}$ and $c_{k}$ induces an isomorphism
$\pi_{i}C_{k}T_n\rtarr \pi_{i}T_n$ for $i>n+k$.
\end{lem}

The following observation is the key to the proofs of Lemmas \ref{lemma1} and \ref{lemma2}. 

\begin{lem}\label{key} Let $T$ be a convergent symmetric spectrum. For each integer $k$ and 
based space $X$, $C_kT\sma X$ is a convergent symmetric spectrum. 
\end{lem}
\begin{proof}
Assume that the structure map $\tilde{\si}: T_n\rtarr \OM T_{n+1}$ is an 
$(n+\la(n))$-equi\-va\-lence, where $\la(n)$ is a nondecreasing 
sequence with infinite limit. The structure maps of $C_kT$ then have the same
property. An easy diagram chase shows that the $n$th structure map of 
$C_kT\sma X$ factors as the composite
$$C_kT_n\sma X \overto{\tilde{\si}\sma\id} (\OM C_kT_{n+1})\sma X
\rtarr  \OM(C_kT_{n+1}\sma X).$$
The first map is an $(n+\lambda(n))$-equivalence and for $n+k>1$ the
second map is a $(2n+2k+1)$-equivalence, by Lemma \ref{easy} below.
The composite is therefore an $(n + \min(n+2k+1,\la(n)))$-equivalence.
\end{proof}

\begin{lem}\label{easy}
If $Y$ is $r$-connected, then the canonical map 
$(\OM Y)\sma X \rtarr \OM (Y\sma X)$ is a $(2r-1)$-equivalence.
\end{lem}
\begin{proof}
This follows from the Freudenthal suspension theorem applied to $\OM Y \sma X$ 
and to $\OM Y$, together with the factorization of the cited map as the composite
$$\OM Y\sma X \overto{\et} \OM\SI(\OM Y\sma X)
\iso \OM(\SI\OM Y\sma X) \overto{\OM(\epz\sma\id)} \OM(Y\sma X),$$
where $\et$ and $\epz$ are the unit and counit of the $(\SI,\OM)$-adjunction.
\end{proof}

We need the following general result about homotopy colimits over $I$ 
to prove Lemmas \ref{lemma1} and \ref{lemma2}.

\begin{prop}\label{connectivity}
Let $F$ and $F'$ be functors $I\rtarr \sT$ and let $\ph: F\rtarr F'$  be a 
natural transformation.  Assume that
$\ph({\mathbf n}):F({\mathbf n}) \rtarr  F'({\mathbf n})$ is a 
$\la(n)$-equivalence, where $\lambda(n) \leq \lambda(n+1)$ and $\lim_n
\la(n)=\infty$. Then the induced map
$$\ph_*: \hocolim_I F \rtarr \hocolim_I F'$$
is a weak equivalence.
\end{prop}

The proof of this proposition is delayed to section~\ref{HOCO} where
we discuss homotopy colimits.

The following simple observation is used several times in the following proofs
of Lemmas~\ref{lemma1} and~\ref{lemma2}. 

\begin{lem}\label{telOmega}
For any sequence of based spaces $Y_n$ and based maps $Y_n\rtarr Y_{n+1}$, 
the evident maps $\OM Y_n\rtarr \OM\tel_n Y_n$ induce a weak equivalence
$$\tel_n \OM Y_n \rtarr \OM \tel_n Y_n.$$
\end{lem}

\begin{proof}[Proofs of Lemmas~\ref{lemma1} and \ref{lemma2}]
We are given a convergent symmetric spectrum $T$.
We first show that both lemmas hold when $T$ is replaced by $C_kT$ 
for any fixed integer $k$. We show that $\bQ C_{k}T$ is a
symmetric $\OM$-spectrum, and the map 
\begin{equation}\label{C1}
C_k T\sma X \rtarr \bQ C_{k}T
\end{equation}
is a $\pi_*$-isomorphism.

The group $\pi_{j}\bQ C_{k}T_n$ 
is the colimit of the sequence
$$\diagram
\pi_{j}(C_{k}T_n\sma X)\rto
& \pi_{j+1}(C_{k}T_{n+1}\sma X)\rto 
& \pi_{j+2}(C_{k}T_{n+2}\sma X)\rto & \cdots.
\enddiagram  $$
We can calculate the maps $\pi_{j}\bQ C_kT_n\rtarr \pi_{j+1}\bQ C_kT_{n+1}$
induced by the structure maps $\bQ C_kT_n\rtarr 
\OM(\bQ C_{k}T_{n+1})$ as maps of colimits
arising from the diagrams
$$ \diagram
\pi_{j}(C_{k}T_n\sma X)\rto\dto
&\pi_{j+1}(C_{k}T_{n+1}\sma X)\rto\dto 
&\cdots\\
\pi_{j+1}(C_{k}T_{n+1}\sma X)\rto
&\pi_{j+2}(C_{k}T_{n+2}\sma X)\rto 
&\cdots.
\enddiagram  $$
Since $C_{k}T\sma X$ is convergent, the vertical maps in each such diagram
are eventually isomorphisms and so induce isomorphisms of the colimits.  Thus each 
structure map $\bQ C_kT_n\rtarr \OM \bQ C_kT_{n+1}$ is
a weak equivalence. The map on homotopy groups induced by (\ref{C1}) is the map
$$\colim_n\pi_{j+n}(C_kT_n\sma X) \rtarr \colim_{m,n}\pi_{j+n+m}(C_kT_{n+m}\sma X)$$
obtained by mapping the terms in the source to the terms with $m=0$ in the target. Extending 
the diagram above vertically and arguing similarly, we see that this map is an isomorphism. 
We should note here that the subtlety mentioned in Warning~\ref{warn} appears
because the vertical maps here differ from the horizontal maps up to 
isomorphism, see~\cite[5.6.3]{hss}.  But this difficulty is avoided here
since $C_kT \sma X$ is convergent. 

For Lemma \ref{lemma2}, since $C_kT \sma X$ is convergent the
colimit, $\colim\pi_{i+k} (C_kT_{n+k}\sma X)$, is attained at 
$\pi_i(C_kT_n \sma X)$ for $i \leq n+ \la(n)$.  
This directly implies the
hypotheses required by Proposition~\ref{connectivity} to show that
\begin{equation}\label{C2}
\hocolim_I \OM^{n}(C_kT_n\sma X)\rtarr
\hocolim_I\OM^{n}\bQ C_kT_n
\end{equation}
is a weak equivalence.

We now deduce that the lemmas hold for $T$. We first show that $\bQ T$ is a symmetric  $\OM$-spectrum. 
Using the maps $C_{-k}T\rtarr C_{-(k+1)}T$, we define
$T' = \tel_k C_{-k}T$.
The compatible maps $C_{-k}T\rtarr T$ induce a level weak equivalence of symmetric spectra
$ c :T'\rtarr T$. 
Observe that $T'\sma X \iso \tel_k(C_{-k}T\sma X)$. Similarly, using the induced maps 
$\bQ C_{-k}T\rtarr \bQ C_{-(k+1)}T,$
we define
$\bQ' T = \tel_{k}\bQ C_{-k}T$.
The compatible maps $\bQ C_{-k}T\rtarr \bQ T$ induce a map of symmetric spectra
$ d :\bQ' T \rtarr \bQ T $
such that the following diagram commutes:
$$ \diagram
T'\sma X\rto \dto_{c\sma\id} & \bQ' T\dto^{d}\\
T\sma X \rto & \bQ T.
\enddiagram $$
The top horizontal arrow is a $\pi_*$-isomorphism since it is obtained by passage
to telescopes from the $\pi_*$-isomorphisms (\ref{C1}). The map $c\sma\id$ is a level weak equivalence since $c$ is a level weak equivalence. Using Lemma \ref{telOmega}, we see that 
$\bQ' T$ is a symmetric  $\OM$-spectrum since each $\bQ C_kT$ is a symmetric  $\OM$-spectrum. 
Thus it suffices to show that $d$ is a level weak equivalence to conclude both that 
$\bQ T$ is a symmetric  $\OM$-spectrum and that the bottom horizontal arrow in the diagram 
is a $\pi_*$-isomorphism, giving Lemma \ref{lemma1}. On passage to telescopes, the
weak equivalences
$$\tel_{k}\OM^{m}(C_{-k}T_{n + m}\sma X)\rtarr 
\OM^{m}\tel_{k}C_{-k}T_{n + m}\sma X $$
given by Lemma \ref{telOmega} induce a weak equivalence
\begin{eqnarray*}
\bQ' T & \iso & \tel_{m}\tel_{k}\OM^{m}(C_{-k}T_{n + m}\sma X)\\
& \longrightarrow &
\tel_{m}\OM^{m}(\tel_{k}C_{-k}T_{n + m}\sma X) = \bQ T'.
\end{eqnarray*}
This is the $n$th term of a level weak equivalence
$ b: \bQ' T\rtarr \bQ T'$.
It is easily seen that $d$ factors as the composite
$$\bQ' T\overto{b} \bQ T'\overto{\bQ c} \bQ T.$$ 
Since $c$ is a level weak equivalence, so is $\bQ c$.
Therefore $d$ is a level weak equivalence.  This completes the proof
of Lemma \ref{lemma1}.

Finally, to complete the proof of Lemma \ref{lemma2}, we observe that the level weak 
equivalences $c$ and $d$ induce the maps $\overline c$ and $\overline d$ displayed in 
the following commutative diagram:
$$ \diagram
\tel_k\hocolim_I \OM^{n}(C_kT_n\sma X)\rto \dto_{\iso}
& \tel_k\hocolim_I\OM^{n}\bQ C_kT_n \dto^{\iso}\\
\hocolim_I\tel_k \OM^{n}(C_kT_n\sma X)\rto \dto_{\htp}
& \hocolim_I\tel_k\OM^{n}\bQ C_kT_n \dto^{\htp}\\
\hocolim_I \OM^nT'_n\sma X)\rto \dto_{\overline{c}}
& \hocolim_I\OM^{n}(\bQ' T_n) \dto^{\overline{d}}\\
\hocolim_{I}\OM^{n}(T_n\sma X)\rto 
&\hocolim_{I}\OM^{n}\bQ T_n.\\
\enddiagram $$
The top horizontal arrow is a weak equivalence because it is the telescope
of the weak equivalences (\ref{C2}). The vertical arrows labelled $\iso$ are
homeomorphisms obtained by commuting $\tel_k$ with $\hocolim_I$. The vertical
arrows labelled $\htp$ are weak equivalences induced by weak equivalences of
Lemma \ref{telOmega}. The maps $\overline c$ and $\overline d$ are weak 
equivalences since $c$ and $d$ are level weak equivalences. Therefore 
the bottom horizontal arrow is a weak equivalence, which is the conclusion of 
Lemma \ref{lemma2}.
\end{proof}

\section{Postnikov towers and connective covers of symmetric spectra}

We construct the strict coverings promised in Lemma~\ref{covering}
from the functorial strict Postnikov towers given by the following lemma. 

\begin{lem}\label{postnikov}
For each integer $k$, there is a functor $P_{k}$ on the category
of symmetric spectra and a natural transformation $\xi_{k}:
\text{Id}\rtarr P_{k}$ such that for every symmetric spectrum $T$ and every
$n$, the map
$$
\xi_{k}(\mathbf n): T_n\rtarr P_{k}T_n
$$
is a $(k+n+1)$-equivalence and $\pi_{i}P_{k}T_n=0$ for $i>k+n$.
Further, there are natural transformations $\chi_{k}:
P_{k}\rtarr P_{k-1}$ such that $\xi_{k-1}=\chi_{k}\circ \xi_{k}$.
\end{lem}

\begin{proof}[Proof of Lemma~\ref{covering}]
We define the symmetric spectrum $C_{k}T$ to be the (levelwise) homotopy 
fiber of the map $\xi_{k}$.  For each $\mathbf n$ in $\SI  $, we define 
$c_k: C_{k}T_n\rtarr T_n$ to be the left vertical arrow in the
following pullback diagram:
$$
\diagram
C_{k}T_n \dto_{c_{k}} \rto & \map_{\sT}(I,P_{k}T_n)\dto\\
T_n\rto_{\xi_{k}(\mathbf n)} &P_{k}T_n.
\enddiagram
$$
It is clear from the
functoriality of our Postnikov sections that $C_{k}$ defines a functor 
on symmetric spectra and that $c_{k}$ is natural.  Since $\xi_{k}(\mathbf n)$ 
is a $(k+n+1)$-equivalence, $C_{k}T$ is strictly $k$-connected.  Since 
$\pi_{i}P_{k}T_n=0$ for $i>k+n$, $c_{k}(\mathbf n)$ induces an isomorphism 
on homotopy groups in degrees $k+n+1$ and above. The maps $e_{k}$ are induced 
by the maps $\chi_{k}$, and the equation $c_{k-1}=c_{k}\circ e_{k}$ follows 
from the corresponding equation $\xi_{k-1}=\chi_{k}\circ \xi_{k}$.
\end{proof}

We will prove Lemma~\ref{postnikov} by a localization technique analogous to one 
that allows a direct construction of functorial Postnikov towers in the category 
of spaces. We will need the left adjoint $F_n$ to the $n$th space functor from 
symmetric spectra to spaces. Thus we have homeomorphisms of hom spaces
$$\map_{\spec}(F_nX,T)\iso \map_{\sT}(X,T_n).$$
As noted in \cite[4.2]{mmss}, $F_n$ is given explicitly by
$$ (F_{n}X)_m= \SI_{m+}\sma_{\SI_{m-n}}X \sma S^{m-n},$$
where we interpret $S^{m-n}$ as $*$ for $m<n$.

We construct $P_{k}$ using the ``small objects argument'' from model category 
theory. For each $n$, define a functor $L_{k,n}$ on the category of
symmetric spectra and a natural transformation $\et: \Id\rtarr L_{k,n}$ as
follows.  For a symmetric spectrum $T$, let $D_{k,n}T$ be the set
of maps of symmetric spectra $F_nS^{k+n+1}\rtarr T$.  Define
$L_{k,n}T$ and $\et$ by the following pushout diagram in the category of
symmetric spectra:
$$ \diagram
\bigvee_{f\in D_{k,n}T} F_nS^{k+n+1}
\rto^(.7){\vee f}\dto_{\vee F_ni}&T\dto^{\et}\\
\bigvee_{f\in D_{k,n}T} F_nCS^{k+n+1}\rto&L_{k,n}T.
\enddiagram $$
Here $i$ denotes the inclusion of the sphere in the cone.

When $m<n$, the space $F_nS^{n+k+1}_m$ is a point and so
$\eta(\mathbf m): T_m\rtarr L_{k,n}T_m$ is an
isomorphism.  When $m\geq n$, the space $F_nS^{n+k+1}_m$ is
a wedge of $S^{m+k+1}$'s, and so the map $\eta(\mathbf m):
T_m\rtarr L_{k,n}T_m$ induces an isomorphism on homotopy
groups in degrees~$0$ to~$m+k$ inclusive.  Moreover, by adjunction,
$\eta(\mathbf n)$ induces the zero map $\pi_{n+k+1}T_n\rtarr \pi_{n+k+1}L_{k,n}T_n$.
Let $N_{k,n}$ be the telescope of iterates of $\eta$:
$$ N_{k,n}= \tel (\Id \overto{\eta} 
L_{k,n}\overto{\eta} 
L_{k,n}\circ L_{k,n}\overto{\eta}
L_{k,n}\circ L_{k,n}\circ L_{k,n} \overto{\eta}\cdots ). $$
Then $\pi_{n+k+1}N_{k,n}T_n=0$ and there is a natural
transformation $\nu_{k,n}: \Id\rtarr N_{k,n}$ such that $\nu_{k,n}(\mathbf
m): T_m\rtarr N_{k,n}T_m$ induces an
isomorphism on homotopy groups in degrees~$0$ to~$m+k$ inclusive.

Let $N_{k}$ be the telescope of the maps $\nu_{k,n}$:
$$ N_{k}= \tel(\Id\overto{\nu_{k,0}}
N_{k,0}\overto{\nu_{k,1}}
N_{k,1}\circ N_{k,0}\overto{\nu_{k,2}}
N_{k,2}\circ N_{k,1}\circ N_{k,0}\overto{\nu_{k,3}}\cdots ). $$
Then $\pi_{n+k+1}N_{k}T_n=0$ for all $n$ and there is a natural
transformation $\nu_{k}: \Id\rtarr N_{k}$ such that $\nu_{k}(\mathbf
n): T_n\rtarr N_{k}T_n$ induces an
isomorphism on homotopy groups in degrees~$0$ to~$n+k$ inclusive for
all $n$.  

Since the map $\nu_{k}$ is the inclusion of the initial object in a
telescope, each map $\nu_{k}(\mathbf n): T_n\rtarr
N_{k}T_n$ is a cofibration.  Let $P_{k}$ be the colimit of
the levelwise cofibrations $\nu_{j}$ for $j>k$:
$$ P_{k}= \colim (N_{k}\overto{\nu_{k+1}}
N_{k+1}\circ N_{k}\overto{\nu_{k+2}}
N_{k+2}\circ N_{k+1}\circ N_{k}\overto{\nu_{k+3}}\cdots ). $$
There is a natural transformation  $\xi_{k}:\Id\rtarr P_{k}$ induced
by $\nu_{k}$.  It follows from the properties of $N_{k}$ that
each space $P_{k}T_n$ has zero homotopy groups above degree
$n+k$ and that each map $\xi_{k}(\mathbf n)$ is a
$(k+n+1)$-equivalence.  There are natural transformations
$\chi_{k}: P_{k}\rtarr P_{k-1}$ induced by the following map of diagrams:
$$ \spreaddiagramcolumns{1pc}\diagram
&N_{k}\rto^(.45){\nu_{k+1}}\dto_{N_{k}\nu_{k-1}}
&N_{k+1}\circ N_{k}\rto^(.65){\nu_{k+2}}\dto_{N_{k+1}N_{k}\nu_{k-1}}&\cdots \\
N_{k-1}\rto_(.45){\nu_{k}}&N_{k}\circ N_{k-1}\rto_(.45){\nu_{k+1}}
&N_{k+1}\circ N_{k}\circ N_{k-1}\rto_(.65){\nu_{k+2}}&\cdots.
\enddiagram $$
The diagram commutes because of the naturality of the horizontal maps.
We see that $\chi_{k}\circ \xi_{k}=\xi_{k-1}$ since the following
diagram commutes:
$$ \diagram
\Id\rto^{\nu_{k}}\dto_{\nu_{k-1}}
&N_{k}\dto^{N_{k}\nu_{k-1}}\\
N_{k-1}\rto_{\nu_{k}}&N_{k}\circ N_{k-1}.
\enddiagram $$
The functors $P_{k}$ and natural transformations $\xi_{k}$ and $\chi_{k}$ have the 
properties specified in Lemma ~\ref{postnikov}.

\section{Homotopy colimits of spaces and symmetric spectra}\label{HOCO}

We first describe what we have used about homotopy colimits of spaces.

Let $\J$ be a small discrete category, thought of as an indexing category.
For a functor $F:\J\rtarr \sT$, the homotopy colimit of $F$ is defined
in terms of the one-sided categorical bar construction $B(*,-,-)$ as
\[\hocolim_\J F = B(*,\J,F).\]
See \cite[X.3]{ekmm} or \cite[\S 12]{class}. The same definition is given
in different language in \cite[XII\S2]{bk}.
This is the geometric realization of a simplicial space, and the inclusion
of its subspace of zero simplices gives a map
\[ \bigvee_{c\in\, \text{Ob}\,\J} F(c)\rtarr \hocolim_\J F.
\]
Homotopy colimits are functorial in $F$: a natural transformation $\ph: F\rtarr F'$ of
functors $\J\rtarr \sT$ induces a map 
\[\ph_*: \hocolim_\J F\rtarr \hocolim_\J F'.\]
The following main property of homotopy colimits, that they are homotopy
invariant, can be proven as in \cite[A.4]{MayPerm}.

\begin{lem}\label{hoho}
If $\ph$ is an levelwise $n$-equivalence, weak equivalence, or homotopy equivalence,
then $\ph_*$ is an $n$-equivalence, weak equivalence, or homotopy equivalence.
\end{lem}

We have used the following easy consequence.

\begin{lem}\label{includezero}
Let $F:\J\rtarr \sT$ be a functor such that $F(\al)$ is a weak equivalence for every 
morphism $\al$ in $\J$ and let $\J$ have an initial object $\ul 0$.  Then the inclusion
$$ F(\ul 0)\rtarr \hocolim_{\J}F $$
is a weak equivalence.
\end{lem}
\begin{proof}
Let $E$ be the constant functor that takes the value $F(\ul 0)$.
By inspection of definitions, $\hocolim_{\J}E$ is homeomorphic to
$F(\ul 0)\sma B\J_+$, and $B\J$ is contractible since $J$ has an
initial object. Therefore the inclusion $F(\ul 0)\rtarr\hocolim_{\J}E$
is a homotopy equivalence. Applying $F$ to the unique map from $\ul 0$ 
to each object of $\J$, we obtain a natural transformation $\ph : E\rtarr F$. 
By the previous lemma, the induced map $ \ph_*: \hocolim_{\J} E\rtarr \hocolim_{\J}F$ 
is a weak equivalence.  Clearly $\ph_*$ restricts to the inclusion on $F(\ul 0)$.
\end{proof}

Homotopy colimits are also functorial in $\J$. Given indexing categories $\J$ and $K$ and
a functor $f: \J\rtarr K$, we obtain an induced map 
\[f_*: \hocolim_\J F\com f \rtarr \hocolim_K F.\]
We have used the following result relating natural transformations to homotopies.

\begin{lem}\label{Nat} 
Let $\et:f\rtarr g$ be a natural transformation between functors
$f,g: \J\rtarr K$. Then the following diagram is naturally homotopy commutative:
$$\diagram
\hocolim_\J F\com f \rrto^{(F\et)_*} \drto_{f_*}& & \hocolim_J F\com g \dlto^{g_*}\\
& \hocolim_K F &
\enddiagram$$
\end{lem}
\begin{proof}
Let $I$ be the category with two objects $[0]$ and $[1]$ and one non-identity
morphism $[0]\rtarr [1]$. Then $BI \iso [0,1]$. It is standard that $\et$ determines
and is determined by a functor $\overline{\et}: \J\times I\rtarr K$ that restricts to $f$ on
$\J\times\{[0]\}$ and to $g$ on $\J\times\{[1]\}$. Let $\pi: \J\times I\rtarr \J$
be the projection. Then $\et$ also determines a natural transformation 
$\tilde\et: f\com \pi\rtarr \overline{\et}$. For an object $c$ of $\J$, 
$\tilde{\et}(c,[0]) = \id: f(c)\rtarr f(c)$ and $\tilde{\et}(c,[1]) = \et: f(c)\rtarr g(c)$.
The required natural homotopy $h: f_*\htp g_*\com(F\et)_*$ is the composite
\begin{eqnarray*}
(\hocolim_\J F\com f)\times BI & \iso & \hocolim_{\J\times I} F\com f\com \pi\\
& \overto{(F\tilde\et)_*} & \hocolim_{\J\times I} F\com \overline{\et}\\
& \overto{\overline{\et}_*} & \hocolim_K F,
\end{eqnarray*}
where the first isomorphism is a direct inspection of definitions.
\end{proof}

We now turn to the 
proof of Proposition~\ref{connectivity} which depends on the expression of 
homotopy colimits over $I$ as homotopy colimits over certain subcategories of 
$I$.
Let $I_{n}$ be the full subcategory of $I$ consisting of the objects 
$\mathbf{m}$ for $m\geq n$ and let $f_{n}$ be the inclusion of $I_{n}$ in $I$. 

\begin{lem}\label{isubn}
Let $F: I\rtarr \sT$ be a functor.  For each $n$, the map
$$ (f_{n})_*: \hocolim_{I_{n}}(F\circ f_{n}) 
\rtarr \hocolim_{I}F $$
induced by $f_n$ is a homotopy equivalence.
\end{lem}
\begin{proof}
Define a functor $g_n: I\rtarr I_n$ by concatenation with
$\mathbf n$; explicitly, $g_n(\mathbf m) = \mathbf{m+n}$ and 
$g_n(i) = i\sqcup\id_{\mathbf n}$ for an injection $i:\mathbf m\rtarr \mathbf m'$. Let
$z_n: \mathbf 0\rtarr \mathbf n$ be the unique map, define 
$\ze(\mathbf m) = \id_{\mathbf m}\sqcup z_n: \mathbf m\rtarr \mathbf{m+n}$, and define
$\et(\mathbf m) = \ze(\mathbf m)$ for $m\geq n$. Then $\ze$ is a natural transformation
$\Id_I\rtarr f_n\com g_n$, $\et$ is a natural transformation $\Id_{I_n}\rtarr g_n\com f_n$,
and $f_n\et = \ze f_n: f_n\rtarr f_n\com g_n\com f_n$. 
Via two applications of Lemma \ref{Nat} and a little diagram chasing, we find 
that the composite
$$\diagram
\hocolim_I F \rto^(0.4){(F\ze)_*} 
& \hocolim_I F\com f_n\com g_n \rto^(0.55){(g_n)_*} & \hocolim_{I_n} F\com f_n
\enddiagram$$
is a homotopy inverse to $(f_n)_*$.
\end{proof}

\begin{proof}[Proof of Proposition \ref{connectivity}]
It suffices to show that $\ph_*$ is an $N$-equivalence for each
$N>0$.  Choose $n$ such that $\lambda(n)>N$.  Then for every object
$\mathbf m$ in $I_{n}$ the map $\ph: F(\mathbf m)\rtarr F'(\mathbf m)$ is an $N$-equivalence. 
By Lemma \ref{hoho}, this implies that
$$ \ph_*: \hocolim_{I_{n}}F\circ f_{n}\rtarr \hocolim_{I_{n}}
F'\circ f_{n} $$
is an $N$-equivalence.  The proposition now follows from Lemma~\ref{isubn}. 
\end{proof}
For a based space $X$, we define $F\sma X$ by 
setting $(F\sma X)(c) = F(c)\sma X$ and have a natural
homeomorphism
\[\hocolim_\J(F\sma X)\iso (\hocolim_\J F)\sma X.\]

We define homotopy colimits of symmetric spectra levelwise. A functor
$F: \J\rtarr \spec$ restricts to give functors $F(\mathbf n): \J\rtarr \sT$ for
each $n\geq 0$, and
\[(\hocolim_\J F)(\mathbf n) = \hocolim_\J F(\mathbf n),\]
and since homotopy colimits commute with the smash product 
$\hocolim_\J F$ is a symmetric spectrum.
All of the statements above remain valid for homotopy colimits of symmetric spectra.

\end{document}